\documentclass[11pt,a4paper]{article}
\usepackage{ifthen,latexsym,amssymb,amsmath,bbm,fixmath}
\usepackage{tikz}

\newcommand{\I}[1]{{\mathbbm #1}}
\renewcommand{\O}[1]{\overline{#1}}

\newcommand{\A}{\mathcal{A}}

\newcommand{\ind}{\mathbbm{1}}

\renewcommand{\mid}{:}
\renewcommand{\ldots}{\hspace{0.9pt}.\hspace{0.3pt}.\hspace{0.3pt}.\hspace{1.5pt}}
\renewcommand{\ge}{\geqslant}
\renewcommand{\le}{\leqslant}

\newif\ifnotesw\noteswtrue
\newcommand{\comment}[1]{\ifnotesw $\blacktriangleright$\ {\sf #1}\ 
  $\blacktriangleleft$ \fi}
\noteswfalse	
\newcommand{\hide}[1]{}

\newcommand{\beq}[1]{\begin{equation}\label{#1}}
\newcommand{\eeq}{\end{equation}}

\newtheorem{theorem}{Theorem}
\newtheorem{lemma}[theorem]{Lemma}
\newtheorem{problem}[theorem]{Problem}

\newtheorem{corollary}[theorem]{Corollary}
\newtheorem{proposition}[theorem]{Proposition}

\newcommand{\bpf}[1][Proof.]{\smallskip\noindent{\it #1} }
\newcommand{\qed}{\nolinebreak\mbox{\hspace{5 true pt}%
  \rule[-0.85 true pt]{3.9 true pt}{8.1 true pt}}}
\newcommand{\epf}{\qed \medskip}

\newtheorem{claim}{Claim}[theorem]

\newcommand{\Z}[1]{\I Z_{#1}} 

\usepackage{scalerel,stackengine}
\stackMath
\newcommand\fourier[1]{%
\savestack{\tmpbox}{\stretchto{%
  \scaleto{%
    \scalerel*[\widthof{\ensuremath{#1}}]{\kern-.6pt\bigwedge\kern-.6pt}%
    {\rule[-\textheight/2]{1ex}{\textheight}}
  }{\textheight}%
}{0.5ex}}%
\stackon[1pt]{#1}{\tmpbox}%
}
\usepackage[T1]{fontenc}
\usepackage[utf8]{inputenc}
\usepackage{authblk}

\title{Minimum number of additive tuples in groups of prime order}

\author[1]{Ostap Chervak}
\author[1]{Oleg Pikhurko\footnote{Supported by ERC
grant~306493 and EPSRC grant~EP/K012045/1.}}
\author[2]{Katherine Staden\footnote{Supported by ERC
grant~306493.}}

\affil[1]{Mathematics Institute and DIMAP\\
University of Warwick\\
Coventry CV4 7AL, UK}

\affil[2]{
Mathematics Institute\\ University of Oxford\\ Oxford OX2 6GG, UK}

\begin{document}

\maketitle 

\begin{abstract} For a prime number $p$ and a sequence of integers $a_0,\dots,a_k\in \{0,1,\dots,p\}$, let $s(a_0,\dots,a_k)$ be the minimum number of $(k+1)$-tuples $(x_0,\dots,x_k)\in A_0\times\dots\times A_k$ with $x_0=x_1+\dots + x_k$, over subsets $A_0,\dots,A_k\subseteq\Z{p}$ of sizes $a_0,\dots,a_k$ respectively. An elegant argument of Lev (independently rediscovered by Samotij and Sudakov) shows that there exists an extremal configuration with all sets $A_i$ being intervals of appropriate length, and that the same conclusion also holds for the related problem, reposed by Bajnok, when $a_0=\dots=a_k=:a$ and  $A_0=\dots=A_k$, provided $k$ is not equal 1 modulo~$p$. By applying basic Fourier analysis, we show for Bajnok's problem that if $p\ge 13$ and $a\in\{3,\dots,p-3\}$ are fixed while $k\equiv 1\pmod p$ tends to infinity, then the extremal configuration alternates between at least two affine non-equivalent sets.
\end{abstract}

\section{Introduction}

Let $\Gamma$ be a given finite Abelian group, with the group operation written additively.

For $A_0,\dots,A_k\subseteq\Gamma$, let $s(A_0,\dots,A_k)$ be the number of $(k+1)$-tuples $(x_0,\dots,x_k)\in A_0\times\dots\times A_k$ with $x_0=x_1+\dots+x_k$. If $A_0=\dots=A_k:=A$, then we use the shorthand $s_k(A):=S(A_0,\dots,A_k)$. For example, $s_2(A)$ is the number of \emph{Schur triples} in $A$, that is, ordered triples $(x_0,x_1,x_2)\in A^3$ with $x_0=x_1+x_2$. 

For integers $n\ge m\ge 0$, let $[m,n]:=\{m,m+1,\dots,n\}$ and 
 $[n]:=[0,n-1]=\{0,\dots,n-1\}$.
For a sequence $a_0,\dots,a_k\in [\,|\Gamma| +1\,] = \lbrace 0,1,\ldots, |\Gamma|\rbrace$, let $s(a_0,\dots,a_k;\Gamma)$ be the minimum of $s(A_0,\dots,A_k)$ over subsets $A_0,\dots,A_k\subseteq \Gamma$ of sizes $a_0,\dots,a_k$ respectively. Additionally, for $a\in [0,p]$, let $s_k(a;\Gamma)$ be the minimum of $s_k(A)$ over all $a$-sets $A\subseteq \Gamma$.

The question of finding the maximal size of a sum-free subset of $\Gamma$ (i.e.\ the maximum $a$ such that $s_2(a;\Gamma)=0$) originated in a paper of Erd\H os~\cite{Erdos65} in 1965 and took 40 years before it was resolved in full generality by Green and Ruzsa~\cite{GreenRuzsa05}. 

In this paper, we are interested in the case where $p$ is a fixed prime and the underlying group $\Gamma$ is taken to be $\Z p$, the cyclic group of order $p$,
which we identify with the additive group of residues modulo $p$ (also using the multiplicative structure on it when this is useful).

Lev~\cite{Lev01duke} solved the problem of finding $s_k(a_0,\ldots,a_k;\Z p)$, where $p$ is prime (in the equivalent guise of considering solutions to $x_1+\ldots + x_k=0$).\footnote{We learned of Lev's work after the publication of this paper. For completeness, we still provide a proof of Theorem~\ref{th:main} in Section~\ref{knot=1}, which is essentially the same as the original proof of Lev's more general result, and which was rediscovered in~\cite{SamotijSudakov16pmcps}.}
For $I\subseteq\Z p$ and $x,y\in\Z p$, write $x\cdot I+y:=\{x\cdot z+y\mid z\in I\}$.

\begin{theorem}\cite{Lev01duke} \label{th:main} For arbitrary $k\ge 1$ and $a_0,\dots,a_k\in [0,p]$, there is $t\in\Z p$ such that 
 $$
 s(a_0,\dots,a_k;\Z p)=s([a_0]+t,[a_1],\dots,[a_k]; \Z p).\qed
 $$
\end{theorem}

Huczynska, Mullen and Yucas~\cite{HuczynskaMullenYucas09jcta} found a new proof of the $s_2(a; \Z p)$-problem, while also addressing some extensions. Samotij and Sudakov~\cite{SamotijSudakov16pmcps} rediscovered Lev's proof of the $s_2(a ; \Z p)$-problem and showed that, when $s_2(a ,\Z p)>0$, then the $a$-sets that achieve the minimum
are exactly those of the form $\xi\cdot I$ with $\xi\in\Z{p}\setminus\{0\}$, where $I$ consists of the residues modulo $p$ of $a$ integers closest to $\frac{p-1}2\in\I Z$. Each such set is an arithmetic progression; its difference can be any non-zero value but the initial element has to be carefully chosen.
(By an \emph{$m$-term arithmetic progression} (or \emph{$m$-AP} for short) we mean a set of the form $\{x,x+d,\dots,x+(m-1)d\}$ for some $x,d\in\Z p$ with $d\not=0$. We call $d$ the \emph{difference}.)
Samotij and Sudakov~\cite{SamotijSudakov16pmcps} also solved the $s_2(a)$-problem for various groups $\Gamma$.
Bajnok~\cite[Problem~G.48]{Bajnok18acmrp} suggested the more general problem of considering $s_k(a;\Gamma)$. 
This is wide open in full generality. 

This paper concentrates on the case $\Gamma = \Z p$, for $p$ prime, and the sets which attain equality in Theorem~\ref{th:main}.
In particular, we write $s(a_0,\dots,a_k):= s(a_0,\dots,a_k;\Z p)$ and $s_k(a):=s_k(a;\Z p)$. Since the case $p=2$ is trivial, let us assume that $p\ge 3$. 
Since  
 \begin{equation}\label{eq:equiv}
 s(A_0,\dots,A_k)=s(\xi\cdot A_0+\eta_0,\dots,\xi\cdot A_k+\eta_k),\quad\mbox{for $\xi\not=0$ and $\eta_0=\eta_1+\dots+\eta_k$},
 \end{equation} 
Theorem~\ref{th:main} shows that, for any difference $d$, there is at least one extremal configuration consisting of $k+1$ arithmetic progressions with the same difference $d$. 

In particular, if $a_0=\dots=a_k=:a$, then one extremal configuration consists of $A_1=\dots=A_k=[a]$ and $A_0=[t,t+a-1]$ for some $t\in\Z p$.
Given this, one can write down some formulas for $s(a_0,\ldots,a_k)$ in terms of $a_0,\ldots,a_k$ involving summation (based on~\eqref{eq:s} or a version of~\eqref{eq:sk(A)}) but there does not seem to be a closed form in general.

 If $k\not\equiv1\pmod p$, then by taking $\xi:=1$, $\eta_1:=\dots:=\eta_k:=-t(k-1)^{-1}$, and $\eta_0:=-kt(k-1)^{-1}$ in~\eqref{eq:equiv}, we can get another extremal configuration where all sets are the same: $A_0+\eta_0=\dots=A_k+\eta_k$. Thus Theorem~\ref{th:main} directly implies the following corollary.

\begin{corollary}\label{cr:main} For every $k\ge 2$ with $k\not\equiv1\pmod p$ and $a\in [0,p]$, there is $t\in\Z p$ such that $s_k(a)=s_k([t,t+a-1])$.\qed\end{corollary}

Unfortunately, if $k\ge 3$, then there may be sets $A$ different from APs that attain equality in Corollary~\ref{cr:main} with $s_k(|A|)>0$ (which is in contrast to the case $k=2$). For example, our (non-exhaustive) search showed that this happens already for $p=17$, when 
$$
 s_3(14)=2255=s_3([-1,12])=s_3([6,18]\cup\{3\}).
 $$
 Also, already the case $k=2$ of the more general Theorem~\ref{th:main} exhibits extra solutions. Of course, by analysing the proof of Theorem~\ref{th:main} or Corollary~\ref{cr:main} one can write a necessary and sufficient condition for the cases of equality. We do this in Section~\ref{knot=1}; in some cases this condition can be simplified.

The first main result of this paper is 
to describe the extremal sets for Corollary~\ref{cr:main} when $k \not\equiv 1 \pmod p$ is sufficiently large.
The proof uses basic Fourier analysis on $\Z p$.

 \begin{theorem}\label{th:knot1} Let a prime $p\ge 7$ and an integer $a\in [3,p-3]$ be fixed, and let $k\not\equiv1\pmod p$ be sufficiently large. Then
there exists $t \in \Z p$ for which the only $s_k(a)$-extremal sets are $\xi\cdot[t,t+a-1]$ for all non-zero $\xi \in \Z p$. 
\end{theorem}

\begin{problem} Find a `good' description of all extremal families for Corollary~\ref{cr:main} (or perhaps Theorem~\ref{th:main}) for $k\ge 3$.\end{problem}

While Corollary~\ref{cr:main} provides an example of an $s_k(a)$-extremal set for $k\not\equiv1\pmod p$, the case $k\equiv1\pmod p$ of the $s_k(a)$-problem turns out to be somewhat special. Here, translating a set $A$ has no effect on the quantity $s_k(A)$. More generally, let $\A$ be the group of all invertible affine transformations of $\Z p$, that is, it consists of maps $x\mapsto \xi\cdot x+\eta$, $x\in\Z p$,  for $\xi,\eta\in \Z p$ with $\xi\not=0$. Then \begin{equation}\label{eq:equiv1}
 s_k(\alpha(A))=s_k(A),\quad \mbox{for every $k\equiv1\!\!\pmod p$\ \ and\  \  $\alpha\in\A$}.
 \end{equation}
 Let us call two subsets $A,B\subseteq \Z p$ \emph{(affine) equivalent} if there is $\alpha\in \A$ with $\alpha(A)=B$. By~\eqref{eq:equiv1}, we need to consider sets only up to this equivalence. 
 Trivially,  any two subsets of $\Z p$ of size $a$ are equivalent if $a \leq 2$ or $a \geq p-2$.

 Our second main result, is to describe the extremal sets when $k \equiv 1 \pmod p$ is sufficiently large, again using Fourier analysis on $\Z p$.

\begin{theorem}\label{th:k1} Let a prime $p\ge 7$ and an integer $a\in [3,p-3]$ be fixed, and let $k\equiv1\pmod p$ be sufficiently large. Then the following statements hold for the $s_k(a)$-problem.
 \begin{enumerate}
 \item If $a$ and $k$ are both even, then $[a]$ is the unique (up to affine equivalence) extremal set.
 \item If at least one of $a$ and $k$ is odd, define $I':=[a-1]\cup\{a\}=\{0,\dots,a-2,a\}$. Then
  \begin{enumerate}
  \item $s_k(a)<s_k([a])$ for all large $k$;
  \item $I'$ is the unique extremal set for infinitely many $k$;
  \item $s_k(a)<s_k(I')$ for infinitely many $k$, provided there are at least three non-equivalent $a$-subsets of $\Z p$.
  \end{enumerate}
 \end{enumerate}

\end{theorem}

It is not hard to see that there are at least three non-equivalent $a$-subsets of $\Z p$ if and only if $p\ge 13$ and $a\in [3,p-3]$, or $p\ge 11$ and $a\in [4,p-4]$. Thus Theorem~\ref{th:k1} characterises pairs $(p,a)$ for which there exists an $a$-subset $A$ which is $s_k(a)$-extremal for \emph{all} large $k\equiv1\pmod p$.

\begin{corollary} Let $p$ be a prime and $a\in[0,p]$. There is an $a$-subset $A\subseteq \Z p$ with $s_k(A)=s_k(a)$ for all large $k\equiv1\pmod p$ if and only if $a\le 2$, or $a\ge p-2$, or $p\in\{7,11\}$ and $a=3$.\qed
\end{corollary}

As is often the case in mathematics, a new result leads to further open problems.

\begin{problem} Given $a\in[3,p-3]$, find a `good' description of all $a$-subsets of $\Z p$ that are $s_k(a)$-extremal for at least one (resp.\ infinitely many) values of $k\equiv1\pmod p$.\end{problem}

\begin{problem} Is it true that for every $a\in[3,p-3]$ there is $k_0$ such that for all $k\ge k_0$ with $k\equiv 1\pmod p$, any two $s_k(a)$-extremal sets are affine equivalent?\end{problem}

\section{Proof of Theorem~\ref{th:main}}\label{knot=1}

For completeness, here we prove Theorem~\ref{th:main}, which is a special case of Theorem~1 in~\cite{Lev01duke}.

Let $A_1,\dots,A_k$ be subsets of $\Z p$. Define $\sigma(x;A_1,\dots,A_k)$ as the number of $k$-tuples $(x_1,\dots,x_k)\in A_1\times\dots\times A_k$ with $x=x_1+\dots+x_k$. Also, for an integer $r\ge 0$, let
 \begin{eqnarray*}
 N_r(A_1,\dots,A_k)&:=&\{x\in\Z p\mid \sigma(x;A_1,\dots,A_k)\ge r\},\\
 n_r(A_1,\dots,A_k)&:=&|N_r(A_1,\dots,A_k)|.
 \end{eqnarray*}

These notions are related to our problem because of the following easy identity:
 \begin{equation}\label{eq:s}
 s(A_0,\dots,A_k)=\sum_{r=1}^\infty |A_0\cap N_r(A_1,\dots,A_k)|.
 \end{equation}

Let an \emph{interval} mean an arithmetic progression with difference $1$, i.e.\ a subset $I$ of $\Z p$ of form $\{x,x+1,\dots,x+y\}$. Its \emph{centre} is $x+y/2\in \I Z_p$; it is unique if $I$ is \emph{proper} (that is, $0<|I|<p$).

Note the following easy properties of the sets $N_r$:
 \begin{enumerate}
 \item These sets are nested: 
  \begin{equation}\label{eq:nested}
  N_0(A_1,\dots,A_k)=\Z p\supseteq N_1(A_1,\dots,A_k)\supseteq N_2(A_1,\dots,A_k)\supseteq \dots
  \end{equation}
 \item If each $A_i$ is an interval with centre $c_i$, then $N_r(A_1,\dots,A_k)$ is an interval with centre $c_1+\dots+c_k$.
  \end{enumerate}

We will also need the following result of Pollard~\cite[Theorem~1]{Pollard75}.

\begin{theorem}\label{th:Pollard} Let $p$ be a prime, $k\ge 1$, and $A_1,\dots,A_k$ be subsets of $\Z{p}$ of sizes $a_1,\dots,a_k$. Then for every integer $r\ge 1$, we have
 $$
 \sum_{i=1}^r n_i(A_1,\dots,A_k)\ge \sum_{i=1}^r n_i([a_1],\dots,[a_k]).\qed
 $$
\end{theorem}

\smallskip

\bpf[Proof of Theorem~\ref{th:main}] Let $A_0,\dots,A_k$ be some extremal sets for the $s(a_0,\dots,a_k)$-problem. We can assume that $0<a_0<p$, because 
$s(A_0,\dots,A_k)$ is $0$ if $a_0=0$ and $\prod_{i=1}^ka_i$ 
if $a_0=p$, regardless of the choice of the sets $A_i$.

Since  $n_0([a_1],\dots,[a_k])=p>p-a_0$ while $n_r([a_1],\dots,[a_k])=0<p-a_0$ when, for example, $r>\prod_{i=1}^{k-1} a_i$, there is a (unique) integer $r_0\ge 0$ such that
 \begin{eqnarray}
 n_{r}([a_1],\dots,[a_k])&>&p-a_0,\quad\mbox{all $r\in [0,r_0]$,}\label{eq:r01}\\
  n_{r}([a_1],\dots,[a_k])&\le &p-a_0,\quad\mbox{all integers $r\ge r_0+1$.}\label{eq:r02}
 \end{eqnarray}
 
The nested intervals $N_1([a_1],\dots,[a_k])\supseteq N_2([a_1],\dots,[a_k])\supseteq\ldots$ have the same centre $c:=((a_1-1)+\dots+(a_k-1))/2$. Thus there is a translation $I:=[a_0]+t$ of $[a_0]$, with $t$ independent of $r$, which has as small as possible intersection with each $N_r$-interval above given their sizes, that is,
 \begin{equation}\label{eq:intersection}
  |I\cap N_r([a_1],\dots,[a_k])|=\max\{\,0,\,n_r([a_1],\dots,[a_k])+a_0-p\,\},\quad \mbox{for all $r\in\I N$}.
  \end{equation}

 This and Pollard's theorem give the following chain of inequalities:
 \begin{eqnarray*}
 s(A_0,\dots,A_k)&\stackrel{\eqref{eq:s}}{=}& \sum_{i=1}^\infty |A_0\cap N_i(A_1,\dots,A_k)|\\
  &\ge & \sum_{i=1}^{r_0} |A_0\cap N_i(A_1,\dots,A_k)|\\
   &\ge & \sum_{i=1}^{r_0} (n_i(A_1,\dots,A_k)+a_0-p)\\
   &\stackrel{\mathrm{Thm~\ref{th:Pollard}}}{\ge} & \sum_{i=1}^{r_0} (n_i([a_1],\dots,[a_k])+a_0-p)\\
   &\stackrel{\eqref{eq:r01}-\eqref{eq:r02}}{=} & \sum_{i=1}^{\infty} \max\{\,0,\, n_i([a_1],\dots,[a_k])+a_0-p\,\}\\
   &\stackrel{\eqref{eq:intersection}}=&   \sum_{i=1}^{\infty}|I\cap N_i([a_1],\dots,[a_k])|\\
   &\stackrel{\eqref{eq:s}}=& s(I,[a_1],\dots,[a_k]),
   \end{eqnarray*}
   giving the required.\epf

Let us write a necessary and sufficient condition for equality in Theorem~\ref{th:main} in the case $a_0,\dots,a_k\in [1,p-1]$. Let $r_0\ge 0$ be defined by \eqref{eq:r01}--\eqref{eq:r02}. Then, by~\eqref{eq:nested}, a sequence $A_0,\dots,A_k\subseteq \Z p$ of sets of sizes respectively $a_0,\dots,a_k$ is extremal if and only if
 \begin{eqnarray}
 A_0\cap N_{r_0+1}(A_1,\dots,A_k)&=&\emptyset,\label{eq:empty}\\
 A_0\cup N_{r_0}(A_1,\dots,A_k)&=& \Z p,\label{eq:whole}\\
 \sum_{i=1}^{r_0} n_i(A_1,\dots,A_k)&=& \sum_{i=1}^{r_0} n_i([a_1],\dots,[a_k]).\label{eq:PollEq}
 \end{eqnarray}

Let us now concentrate on the case $k=2$, trying to simplify the above condition. We can assume that no $a_i$ is equal to 0 or $p$ (otherwise the choice of the other two sets has no effect on $s(A_0,A_1,A_2)$ and every triple of sets of sizes $a_0$, $a_1$ and $a_2$ is extremal). Also, as in~\cite{SamotijSudakov16pmcps}, let us exclude the case $s(a_0,a_1,a_2)=0$, as then there are in general many extremal configurations. Note that $s(a_0,a_1,a_2)=0$ if and only if $r_0=0$; also, 
by the Cauchy-Davenport theorem (the special case $k=2$ and $r=1$ of Theorem~\ref{th:Pollard}), this is equivalent to $a_1+a_2-1\le p-a_0$. 
Assume by symmetry that $a_1\le a_2$. 
Note that~\eqref{eq:r01} implies that $r_0\le a_1$.

The condition in~\eqref{eq:PollEq} states that we have equality in Pollard's theorem. A result of Nazarewicz, O'Brien, O'Neill and Staples~\cite[Theorem~3]{NazarewiczObrienOneillStaples07} characterises when this happens (for $k=2$), which in our notation is the following. 
 
 \begin{theorem}\label{th:NazarewiczObrienOneillStaples07}  For $k=2$ and $1\le r_0\le a_1\le a_2<p$, we have equality in~\eqref{eq:PollEq} if and only if at least one of the following conditions holds:
   \begin{enumerate}
   \item\label{it:1} $r_0=a_1$,
   \item\label{it:2} $a_1+a_2\ge p+r_0$,
   \item\label{it:3} $a_1=a_2=r_0+1$ and $A_2=g-A_1$ for some $g\in \Z{p}$,
   \item\label{it:4} $A_1$ and $A_2$ are arithmetic progressions with the same difference.
 \end{enumerate}
 \end{theorem}

Let us try to write more explicitly each of these four cases, when combined with~\eqref{eq:empty} and~\eqref{eq:whole}. 

First, consider the case $r_0=a_1$. We have $N_{a_1}([a_1],[a_2])=[a_1-1,a_2-1]$ and thus $n_{a_1}([a_1],[a_2])=a_2-a_1+1>p-a_0$, that is, $a_2-a_1\ge p-a_0$. The condition~\eqref{eq:empty} holds automatically since $N_i(A_1,A_2)=\emptyset$ whenever $i>|A_1|$. The other condition~\eqref{eq:whole} may be satisfied even when none of the sets $A_i$ is an arithmetic progression (for example, take $p=13$, $A_1=\{0,1,3\}$, $A_2=\{0,2,3,5,6,7,9,10\}$ and let $A_0$ be the complement of $N_3(A_1,A_2)=\{3,6,10\}$). We do not see any better characterisation here, apart from stating that~\eqref{eq:whole} holds.

Next, suppose that $a_1+a_2\ge p+r_0$. Then, for any two sets $A_1$ and $A_2$ of sizes $a_1$ and $a_2$, we have $N_{r_0}(A_1,A_2)=\Z p$; thus~\eqref{eq:whole} holds automatically. Similarly to the previous case, there does not seem to be a nice characterisation of~\eqref{eq:empty}. For example,~\eqref{eq:empty} may hold even when none of the sets $A_i$ is an AP: e.g.\ let $p=11$, $A_1=A_2=\{0,1,2,3,4,5,7\}$, and let $A_0=\{0,2,10\}$ be the complement of $N_4(A_1,A_2)=\{1,3,4,5,6,7,8,9\}$ (here $r_0=3$).%
 \comment{Let us verify that indeed $r_0=3$. Indeed, $n_3([7],[7])=11$ by above while $N_4([7],[7])=[3,9]$ has $7\le 11-a_0=8$ elements.
 }

Next, suppose that we are in the third case. The primality of $p$ implies that $g\in\Z p$ satisfying $A_2=g-A_1$ is unique and thus $N_{r_0+1}(A_1,A_2)=\{g\}$. Therefore~\eqref{eq:empty} is equivalent to $A_0\not\ni g$. Also, note that if $I_1$ and $I_2$ are intervals of size $r_0+1$, then $n_{r_0}(I_1,I_2)=3$. By the definition of $r_0$, we have
$p-2\le a_0\le p-1$. 
  Thus we can choose any integer $r_0\in [1,p-2]$ and $(r_0+1)$-sets $A_2=g-A_1$, and then let $A_0$ be obtained from $\Z p$ by removing $g$ and at most one further element of $N_{r_0}(A_1,A_2)$. Here, $A_0$  is always an AP (as a subset of $\Z p$ of size $a_0\ge p-2$) but $A_1$ and $A_2$ need not be.
 
Finally, let us show that if $A_1$ and $A_2$ are arithmetic progressions with the same difference $d$ and we are not in Case~1 nor~2 of Theorem~\ref{th:NazarewiczObrienOneillStaples07}, then $A_0$ is also an arithmetic progression whose difference is~$d$. By~\eqref{eq:equiv}, it is enough to prove this when $A_1=[a_1]$ and $A_2=[a_2]$ (and $d=1$). 
Since $a_1+a_2\le p-1+r_0$ and  $r_0+1\le a_1\le a_2$, we have that
 \begin{eqnarray*} N_{r_0}(A_1,A_2)&=&[r_0-1,a_1+a_2-r_0-1]\\
 N_{r_0+1}(A_1,A_2)&=& [r_0,a_1+a_2-r_0-2]
 \end{eqnarray*} 
 have sizes respectively $a_1+a_2-2r_0+1<p$ and $a_1+a_2-2r_0-1>0$. We see that $N_{r_0+1}(A_1,A_2)$ is obtained from the proper interval $N_{r_0}(A_1,A_2)$ by removing its two endpoints. Thus $A_0$, which is sandwiched between the complements of these two intervals by~\eqref{eq:empty}--\eqref{eq:whole}, must be an interval too. (And, conversely, every such triple of intervals is extremal.)

\section{The proof of Theorems~\ref{th:knot1} and~\ref{th:k1}}

 Let us recall the basic definitions and facts of Fourier analysis on $\Z p$. For a more detailed treatment of this case, see e.g.~\cite[Chapter~2]{Terras99faofg}.
Write $\omega := e^{2\pi i /p}$ for the $p^{\mathrm{th}}$ root of unity.
Given a function $f : \Z p \rightarrow \mathbb{C}$, we define its \emph{Fourier transform} to be the function $\fourier{f}:\Z p\to \mathbb{C}$ given by
$$
\fourier{f}(\gamma) := \sum_{x=0}^{p-1} f(x)\, \omega^{-x\gamma}, \qquad\text{for } \gamma \in \Z p.
$$
Parseval's identity states that
 \begin{equation}\label{eq:Parseval}
\sum_{x=0}^{p-1} f(x)\,\overline{g(x)} = \frac{1}{p}\sum_{\gamma=0}^{p-1} \fourier{f}(\gamma)\,\overline{\fourier{g}(\gamma)}.
\end{equation}
The \emph{convolution} of two functions $f,g : \Z p \rightarrow \mathbb{C}$ is given by
$$
(f * g)(x) := \sum_{y=0}^{p-1}f(y)\,g(x-y).
$$
It is not hard to show that the Fourier transform of a convolution equals the product of Fourier transforms, i.e.
\begin{equation}\label{convolution}
\fourier{f_1 * \ldots * f_k} = \fourier{f_1} \cdot\ldots \cdot \fourier{f_k}.
\end{equation}
We write $f^{*k}$ for the convolution of $f$ with itself $k$ times. (So, for example, $f^{* 2} = f * f$.)
Denote by $\ind_A$ the \emph{indicator function} of $A \subseteq \Z p$ which assumes value 1 on $A$ and $0$ on $\Z p\setminus A$. We will call $\fourier{{\ind}}_A(0)=|A|$ the \emph{trivial Fourier coefficient of $A$}. 
Since the Fourier transform behaves very nicely with respect to convolution, it is not surprising that our parameter of interest, $s_k(A)$, can be written as a simple function of the Fourier coefficients of $\ind_A$.
Indeed,
let $ A \subseteq \mathbb{Z}_p$ and $x \in \Z p$.
Then the number of tuples $(a_1,\ldots,a_k) \in A^k$ such that $a_1+\ldots+a_k=x$ (which is $\sigma(x;A,\ldots,A)$ in the notation of Section~\ref{knot=1}) is precisely $\ind^{* k}_A(x)$. The function $s_k(A)$ counts such a tuple if and only if its sum $x$ also lies in $A$.
Thus, 
\begin{equation}\label{eq:sk(A)}
s_k(A) = \sum_{x = 0}^{p-1}\ind^{* k}_A(x)\, \ind_A(x) \stackrel{(\ref{eq:Parseval})}= \frac{1}{p}\sum_{\gamma=0}^{p-1}\fourier{\ind^{* k}_A}(\gamma)\,\overline{\fourier{\ind_A}(\gamma)} \stackrel{(\ref{convolution})}{=} \frac{1}{p} \sum_{\gamma=0}^{p-1} \left(\fourier{\ind_A}(\gamma)\right)^k\, \overline{\fourier{\ind_A}(\gamma)}.
\end{equation}
Since every set $A \subseteq \Z p$ of size $a$ has the same trivial Fourier coefficient (namely $\fourier{\ind_A}(0)=a$), let us re-write~\eqref{eq:sk(A)} as
\beq{eq:SkF1}
p s_k(A)-a^{k+1}= \sum_{\gamma=1}^{p-1} (\fourier{\ind_A}(\gamma))^k\, 
 \O{\fourier{\ind_A}(\gamma)} =: F(A).
 \eeq
Thus we need to minimise $F(A)$  (which is a real number for any $A$) over $a$-subsets $A\subseteq\Z p$. 
To do this when $k$ is sufficiently large, we will consider the largest in absolute value non-trivial Fourier coefficient $\fourier{\ind_{A}}(\gamma)$ of an $a$-subset $A$. Indeed, the term $(\fourier{\ind_A}(\gamma))^k\overline{\fourier{\ind_A}(\gamma)}$ will dominate $F(A)$, so if it has strictly negative real part, then $F(A)<F(B)$ for all $a$-subsets~$B\subseteq \Z p$ with $\max_{\delta\not=0}|\fourier{\ind_B}(\delta)|<|\fourier{\ind_A}(\gamma)|$.

Given $a \in [p-1]$, let
$$
I := [a]=\lbrace 0,\ldots,a-1\rbrace\quad\text{and}\quad I' := [a-1]\cup\lbrace a \rbrace = \lbrace a,\ldots,a-2,a\rbrace.
$$
In order to prove Theorems~\ref{th:knot1} and~\ref{th:k1}, we will make some preliminary observations about these special sets.
The set of $a$-subsets which are affine equivalent to $I$ is precisely the set of $a$-APs.

Next we will show that
\begin{equation}\label{skI}
F(I) = 2\sum_{\gamma=1}^{(p-1)/2} (-1)^{\gamma(a-1)(k-1)} \left|\fourier{\ind_I}(\gamma)\right|^{k+1}\quad\text{if }k \equiv 1 \pmod p.
\end{equation}
Note that $(-1)^{\gamma(a-1)(k-1)}$ equals $(-1)^\gamma$ if both $a,k$ are even and 1 otherwise.
To see~\eqref{skI}, let $\gamma \in \lbrace 1,\ldots,\frac{p-1}{2}\rbrace$ and write $\fourier{\ind_I}(\gamma) = re^{\theta i}$ for some $r >0$ and $0 \leq \theta < 2\pi$.
Then $\theta$ is the midpoint of $0,-2\pi \gamma/p,\ldots, -2(a-1)\gamma\pi/p$,~i.e.
$
\theta =  -\pi(a-1)\gamma/p
$.
Choose $s \in \mathbb{N}$ such that $k = sp+1$.
Then
\begin{equation}\label{termx}
(\fourier{\ind_I}(\gamma))^k \overline{\fourier{\ind_I}(\gamma)} = \left(r e^{-\pi i (a-1)\gamma/p}\right)^k r e^{\pi i (a-1)\gamma/p} = r^{k+1} e^{-\pi i(a-1)\gamma s},
\end{equation}
and $e^{-\pi i(a-1)s}$ equals $1$ if $(a-1)s$ is even, and $-1$ if $(a-1)s$ is odd.
Note that, since $p$ is an odd prime, $(a-1)s$ is odd if and only if $a$ and $k$ are both even.
So~(\ref{termx}) is real, and the fact that $\fourier{\ind_I}(p-\gamma) = \overline{\fourier{\ind_I}(\gamma)}$ implies that the corresponding term for $p-\gamma$ is the same as for $\gamma$.
This gives~\eqref{skI}.
A very similar calculation to~(\ref{termx}) shows that
\begin{equation}\label{skI2}
F(I+t) = \sum_{\gamma=1}^{p-1} e^{-\pi i (2t+a-1)(k-1)\gamma/p}|\fourier{\ind_{I+t}}(\gamma)|^{k+1}\quad\text{for all }k \geq 3.
\end{equation}

Given $r>0$ and $0 \leq \theta < 2\pi$, we write $\arg(re^{\theta i}) := \theta$.

\begin{proposition}\label{angleprop}
Suppose that $p \geq 7$ is prime and $a \in [3,p-3]$.
Then 
$\arg\left(\fourier{\ind_{I'}}(1)\right)$ is not an integer multiple of $\pi/p$.
\end{proposition}

\bpf
Since $\fourier{\ind_A}(\gamma)=-\fourier{\ind_{\Z p\setminus A}}(\gamma)$ for all $A \subseteq \Z p$ and non-zero $\gamma \in \Z p$, we may assume without loss of generality that $a \leq p-a$. Since $p$ is odd, we have $a \leq (p-1)/2$.

Suppose first that $a$ is odd.
Let $m := (a-1)/2$. Then $m \in [1,\frac{p-3}{4}]$.
Observe that translating any $A \subseteq \Z p$ changes the arguments of its Fourier coefficients by an integer multiple of $2\pi/p$.
So, for convenience of angle calculations, here we may redefine $I := [-m,m]$ and $I' := \lbrace -m-1\rbrace\cup[-m+1,m]$.
Also let $I^- := [-m+1,m-1]$, which is non-empty.
The argument of $\fourier{\ind_{I^-}}(1)$ is $0$.
Further, $\fourier{\ind_{I'}}(1) = \fourier{\ind_{I^-}}(1) + \omega^{m+1}+\omega^{-m}$. 
Since $\omega^{m+1},\omega^{-m}$ lie on the unit circle, the argument of $\omega^{m+1}+\omega^{-m}$ is either $\pi/p$ or $\pi+\pi/p$. But the bounds on $m$ imply that it has positive real part, so $\arg(\omega^{m+1}+\omega^{-m})=\pi/p$.
By looking at the non-degenerate parallelogram in the complex plane with vertices $0,\fourier{\ind_{I^-}}(1),\omega^{m+1}+\omega^{-m},\fourier{\ind_{I'}}(1)$, we see that the argument of $\fourier{\ind_{I'}}(1)$ lies strictly between that of $\fourier{\ind_{I^-}}(1)$ and $\omega^{m+1}+\omega^{-m}$, i.e.~strictly between $0$ and $\pi/p$, giving the required.

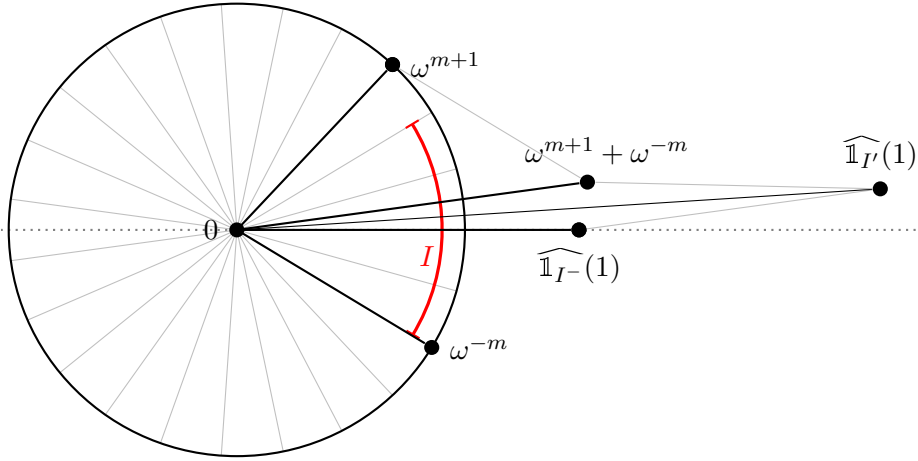
\begin{figure}[h]\label{figure}
\centering

	\begin{tikzpicture}

\clip (-3.1,-3.1) rectangle  (9,3.1);

\draw[thick,gray,dotted] (-3.5,0) -- (9,0);

\foreach \x in {0,15.65217,31.30435,...,350} 
{
\draw[gray!50] (0,0) -- (\x:3);
}

\draw[black,thick] (0,0) circle (3cm); 

\draw (0,0)  node[circle,inner sep=2,fill = black,label=left:{$0$}]  (0) {};
\draw (46.9563:3)  node[circle,inner sep=2,fill = black,label=right:{$\omega^{m+1}$}]  (m+1) {};
\draw (-31.3043:3)  node[circle,inner sep=2, fill=black,label=right:{$\omega^{-m}$}]  (-m) {};


\draw[red,very thick] (31.3043:2.7) arc (31.3043:-31.3043:2.7);
\draw[] (2.5,-0.75)  node[draw=none,label=above:{\textcolor{red}{$I$}}]  () {};
\draw[red,thick] (31.3043:2.6) -- (31.3043:2.8);
\draw[red,thick] (-30.9:2.6) -- (-30.9:2.8);

\begin{scope}[shift=(m+1)]

\draw (-31.3043:3)  node[circle,inner sep=2, fill=black,label=above:{~~~~$\omega^{m+1}+\omega^{-m}$}]  (pip) {};
\draw[gray!50] (0,0) -- (pip);
\end{scope}

\draw (46.9563:3)  node[circle,inner sep=2,fill = black]  (m+1) {};


\draw[black,thick] (0,0) -- (m+1);
\draw[black,thick] (0,0) -- (-m);
\draw[black,thick] (0,0) -- (pip);

\draw (0:4.5)  node[circle,inner sep=2, fill=black,label=below:{$\fourier{\ind_{I^-}}(1)$}]  (I-) {};
\draw[thick] (0,0) -- (I-);

\begin{scope}[shift=(I-)]

\draw[] (7.8261:4)  node[circle,inner sep=2, fill=black,label=above:{$\fourier{\ind_{I'}}(1)$}]  (pipshift) {};
\draw[gray!50] (0,0) -- (pipshift);
\end{scope}

\draw (0,0) -- (pipshift);

\draw[gray!50] (pipshift) -- (pip);
\draw[] (I-)  node[circle,inner sep=2, fill=black]  () {};

	\end{tikzpicture}
\end{figure}

Suppose now that $a$ is even and let $m := (a-2)/2 \in [1,\frac{p-5}{4}]$.
Again without loss of generality we may redefine $I := [-m,m+1]$ and $I' := \lbrace -m-1\rbrace \cup [ -m+1,m+1]$.
Let also $I^- := [-m+1,m]$, which is non-empty.
The argument of $\fourier{\ind_{I^-}}(1)$ is $-\pi/p$.
Further, $\fourier{\ind_{I'}}(1) = \fourier{\ind_{I^-}}(1) + \omega^{m+1}+\omega^{-(m+1)}$.
The argument of $\omega^{m+1}+\omega^{-(m+1)}$ is $0$, so as before the argument of $\fourier{\ind_{I'}}(1)$ is strictly between $-\pi/p$ and $0$, as required.
\epf

We say that an $a$-subset $A$ is a \emph{punctured interval} if $A=I'+t$ or $A = -I'+t$ for some $t \in \Z p$.
That is, $A$ can be obtained from an interval of length $a+1$ by removing a penultimate point.

\begin{lemma}\label{int-equivalence}
Let $p \geq 7$ be prime and let $a \in \lbrace 3,\ldots,p-3\rbrace$.
Then the sets $I,I'\subseteq \Z p$ are not affine equivalent. 
Thus no punctured interval is affine equivalent to an interval.
\end{lemma}
\bpf Suppose on the contrary that there is $\alpha \in \mathcal{A}$ with $\alpha(I')=I$. Let a \emph{reflection} mean an affine map $R_c$ with $c\in\Z p$ that maps $x$ to $-x+c$. 
Clearly, $I=[a]$ is invariant under the reflection $R:=R_{a-1}$. Thus $I'$ is invariant under the map $R':=\alpha^{-1}\circ R\circ \alpha$. As is easy to see, $R'$ is also some reflection and thus preserves the cyclic distances in $\Z p$. So $R'$ has to fix $a$, the unique element of $I'$ with both distance-1 neighbours lying outside of $I'$. Furthermore, $R'$ has to fix $a-2$, the unique element of $I'$ at distance 2 from $a$. However, no reflection can fix two distinct elements  of $\Z p$, a contradiction.
\epf

We remark that the previous lemma can also be deduced from Proposition~\ref{angleprop}. Indeed, for any $A \subseteq \Z p$, the multiset of Fourier coefficients of $A$ is the same as that of $x\cdot A$ for $x \in \Z p\setminus\lbrace 0 \rbrace$, and translating a subset changes the argument of Fourier coefficients by an integer multiple of $2\pi/p$. Thus for every subset which is affine equivalent to $I$, the argument of each of its Fourier coefficients is an integer multiple of $\pi/p$.

Let
$$
\rho(A) := \max_{\gamma \in \Z p \setminus \lbrace 0 \rbrace}|\fourier{\ind}_A(\gamma)|\quad\text{and}\quad R(a) := \left\lbrace \rho(A) : A \in \binom{\Z p}{a}\right\rbrace = \lbrace m_1(a) > m_2(a) > \ldots \rbrace.
$$
Given $j \geq 1$, we say that $A$ \emph{attains} $m_j(a)$, and specifically that \emph{$A$ attains $m_j(a)$ at $\gamma$} if $m_j(a) = \rho(A)=|\fourier{\ind_A}(\gamma)|$.
Notice that, since $\fourier{\ind_A}(-\gamma)=\overline{\fourier{\ind_A}(\gamma)}$, the set $A$ attains $m_j(a)$ at $\gamma$ if and only if $A$ attains $m_j(a)$ at $-\gamma$ (and $\gamma,-\gamma \neq 0$ are distinct values).

As we show in the next lemma, the $a$-subsets which attain $m_1(a)$ are precisely the affine images of $I$ (i.e.~arithmetic progressions), and the $a$-subsets which attain $m_2(a)$ are the affine images of the punctured interval~$I'$.

\begin{lemma}\label{lm:MaxNontriv} Let $p\geq 7$ be prime and let $a\in [3,p-3]$. Then $|R(a)| \geq 2$ and 
\begin{itemize}
\item[(i)] $A \in \binom{\Z p}{a}$ attains $m_1(a)$ if and only if $A$ is affine equivalent to $I$, and every interval attains $m_1(a)$ at $1$ and $-1$ only;
\item[(ii)] $B\in\binom{\Z p}{a}$ attains $m_2(a)$ if and only if $B$ is affine equivalent to $I'$, and every punctured interval attains $m_2(a)$ at $1$ and $-1$ only.
\end{itemize}
\end{lemma}

\bpf
Given $D \in \binom{\Z p}{a}$, we claim that there is some $D_{\rm{pri}} \in \binom{\Z p}{a}$ with the following properties:
\begin{itemize}
\item $D_{\rm{pri}}$ is affine equivalent to $D$;
\item $\rho(D) = |\fourier{\ind_{D_{\rm{pri}}}}(1)|$; and
\item $-\pi/p < \arg\left(\fourier{\ind_{D_{\rm{pri}}}}(1)\right) \leq \pi/p$.
\end{itemize}
Call such a $D_{\rm{pri}}$ a \emph{primary image} of $D$.
Indeed, suppose that $\rho(D) = |\fourier{\ind_D}(\gamma)|$ for some non-zero $\gamma \in \Z p$, and let $\fourier{\ind_D}(\gamma) = r'e^{\theta' i}$ for some $r' > 0$ and $0 \leq \theta' < 2\pi$.
(Note that we have $r'>0$ since $p$ is prime.)
Choose $\ell \in \lbrace 0,\ldots,p-1\rbrace$ and $-\pi/p < \phi \leq \pi/p$ such that $\theta' = 2\pi \ell/p + \phi$. 
Let $D_{\rm{pri}} := \gamma\cdot D + \ell$. Then
$$
|\fourier{\ind_{D_{\rm{pri}}}}(1)| = \left|\sum_{x \in D}\omega^{-\gamma x - \ell}\right| = |\omega^{-\ell} \fourier{\ind_D}(\gamma)| = |\fourier{\ind_D}(\gamma)| = \rho(D),
$$
and
$$
\arg\left(\fourier{\ind_{D_{\rm{pri}}}}(1)\right) = \arg(e^{\theta' i}\omega^{-\ell}) = 2\pi \ell/p + \phi - 2\pi \ell/p = \phi,
$$
as required.

Let $D \subseteq \Z p$ have size $a$ and write $\fourier{\ind_D}(1) = re^{\theta i}$.
Assume by the above that $-\pi/p < \theta \leq \pi/p$.
For all $j \in \Z p$, let
$$
h(j) := \Re(\omega^{-j}e^{-\theta i}) = \cos\left(\frac{2\pi j}{p}+\theta\right),
$$
 where $\Re(z)$ denotes the real part of $z\in\mathbb{C}$.
Given any $a$-subset $E$ of $\Z p$, we have 
\begin{equation}\label{hmbound}
H_D(E) := \sum_{j \in E}h(j) = \Re\left(e^{-\theta i}\sum_{j \in E}\omega^{-j}\right) = \Re\left(e^{-\theta i} \fourier{\ind_E}(1)\right) \leq |\fourier{\ind_E}(1)|.
\end{equation}
Then
\begin{equation}\label{HAA}
H_D(D) = \sum_{j \in D}h(j) = \Re(e^{-\theta i} \fourier{\ind_D}(1)) = r = |\fourier{\ind_D}(1)|.
\end{equation}

Note that $H_D(E)$ is the (signed) length of the orthogonal projection of $\fourier{\ind_E}(1)\in\mathbb{C}$ on the 1-dimensional 
line $\{xe^{i\theta}\mid x\in\I R\}$. As stated in~\eqref{hmbound} and~\eqref{HAA}, $H_D(E)\le |\fourier{\ind_E}(1)|$
and this is equality for $E=D$. (Both of these facts are geometrically obvious.) If $|\fourier{\ind_D}(1)|=m_1(a)$ is maximum, then no 
$H_D(E)$ for an $a$-set $E$ can exceed $m_1(a)=H_D(D)$. Informally speaking, the main idea of the proof is that if we fix the direction $e^{i\theta}$, then the projection length is maximised if we take $a$ distinct elements $j\in \I Z_p$ with the $a$ largest values of $h(j)$, that is, if we take some interval (with the runner-up being a punctured interval). 

Let us provide a formal statement and proof of this now.

\begin{claim}\label{claim}
Let $\mathcal{I}_a$ be the set of length-$a$ intervals in $\Z p$.
\begin{itemize}
\item[(i)]
Let $M_1(D) \subseteq \binom{\Z p}{a}$ consist of $a$-sets $E\subseteq \Z p$  such that $H_D(E) \geq H_D(C)$ for all $C \in \binom{\Z p}{a}$.
Then $M_1(D) \subseteq \mathcal{I}_a$.
\item[(ii)] Let $M_2(D) \subseteq \binom{\Z p}{a}$ be the set of $E \notin \mathcal{I}_a$ for which $H_D(E) \geq H_D(C)$ for all $C \in \binom{\Z p}{a} \setminus \mathcal{I}_a$.
Then every $E \in M_2(A)$ is a punctured interval.
\end{itemize}
\end{claim}

\bpf
Suppose that $0 < \theta < \pi/p$.
Then $h(0) > h(1) > h(-1) > h(2) > h(-2) > \ldots > h(\frac{p-1}{2}) > h(-\frac{p-1}{2})$.
In other words, $h(j_\ell) > h(j_k)$ if and only if $\ell < k$, where $j_m := (-1)^{m-1}\lceil m/2\rceil$.
Letting $J_{a-1} := \lbrace j_0,\ldots,j_{a-2}\rbrace$, we see that
$$
H_D(J_{a-1} \cup \lbrace j_{a-1}\rbrace) > H_D(J_{a-1} \cup \lbrace j_{a}\rbrace) > H_D(J_{a-1} \cup \lbrace j_{a+1}\rbrace),  H_D(J_{a-2} \cup \lbrace j_{a-1},j_a\rbrace) > H_D(J)
$$
for all other $a$-subsets $J$.
But $J_{a-1} \cup \lbrace j_{a-1}\rbrace$ and $J_{a-1} \cup \lbrace j_a\rbrace$ are both intervals, and $J_{a-1} \cup \lbrace j_{a+1}\rbrace$ and $J_{a-2} \cup \lbrace j_{a-1},j_a\rbrace$ are both punctured intervals.
So in this case $M_1(D) := \lbrace J_{a-1}\cup\lbrace j_{a-1}\rbrace\rbrace$ and $M_2(D) \subseteq \lbrace J_{a-1}\cup\lbrace j_{a+1}\rbrace, J_{a-2} \cup \lbrace j_{a-1},j_a\rbrace\rbrace$, as required.

The case when $-\pi/p < \theta < 0$ is almost identical except now $j_\ell := (-1)^\ell\lceil \ell/2\rceil$ for all $0 \leq \ell \leq p-1$.
If $\theta=0$ then $h(0) > h(1) = h(-1) > h(2) = h(-2) > \ldots > h(\frac{p-1}{2}) = h(-\frac{p-1}{2})$.
If $\theta=-\pi/p$ then $h(0)=h(-1) > h(1)=h(-2) > \ldots = h(-\frac{p-1}{2}) > h(\frac{p-1}{2})$.
\epf

\medskip
\noindent
We can now prove part~(i) of the lemma.
Suppose $A \in \binom{\Z p}{a}$ attains $m_1(a)$ at $\gamma \in \Z p \setminus \lbrace 0 \rbrace$.
Then the primary image $D$ of $A$ satisfies $|\fourier{\ind_D}(1)|=m_1(a)=|\fourier{\ind_A}(\gamma)|$.
So, for any $E \in M_1(D)$,
$$
|\fourier{\ind_A}(\gamma)| = |\fourier{\ind_D}(1)| \stackrel{(\ref{HAA})}{=} H_D(D) \leq H_D(E) \stackrel{(\ref{hmbound})}{\leq} |\fourier{\ind_E}(1)|,
$$
with equality in the first inequality if and only if $D \in M_1(D)$.
Thus, by Claim~\ref{claim}(i), $D$ is an interval, and so $A$ is affine equivalent to an interval, as required.
Further, if $A$ is an interval then $D$ is an interval if and only if $\gamma=\pm 1$.
This completes the proof of (i). 

\medskip
\noindent
For~(ii), note that $m_2(a)$ exists since by Lemma~\ref{int-equivalence}, there is a subset (namely $I'$) which is not affine equivalent to $I$. By~(i), it does not attain $m_1(a)$, so $\rho(I') \leq m_2(a)$.
Suppose now that $B$ is an $a$-subset of $\Z p$ which attains $m_2(a)$ at $\gamma \in \Z p \setminus \lbrace 0 \rbrace$.
Let $D$ be the primary image of $B$.
Then $D$ is not an interval.
This together with Claim~\ref{claim}(i) implies that $H_D(D) < H_D(E)$ for any $E \in M_1(D)$.
Thus, for any $C \in M_2(D)$, we have
$$
m_2(a) = |\fourier{\ind_B}(\gamma)| = |\fourier{\ind_D}(1)|  = H_D(D) \leq H_D(C) \leq |\fourier{\ind_C}(1)|.
$$
with equality in the first inequality if and only if $D \in M_2(D)$.
Since $C$ is a punctured interval, it is not affine equivalent to an interval. So the first part of the lemma implies that $|\fourier{\ind_C}(1)|\le m_2(a)$.
Thus we have equality everywhere and so $D \in M_2(D)$.
Therefore $B$ is the affine image of a punctured interval, as required.
Further, if $B$ is a punctured interval, then $D$ is a punctured interval if and only if $\gamma=\pm 1$.
This completes the proof of (ii).
\epf

We will now prove Theorem~\ref{th:knot1}.

\bpf[Proof of Theorem~\ref{th:knot1}.] 
Recall that $p\ge 7$, $a\in [3,p-3]$ and  $k > k_0(a,p)$ is sufficiently large with $k\not\equiv 1\pmod p$. Let $I = [a]$.
Given $t \in \Z p$, write $\rho_t := (\fourier{\ind_{I+t}}(1))^k\overline{\fourier{\ind_{I+t}}(1)}$ as $r_te^{\theta_t i}$, where $\theta_t \in [0,2\pi)$ and $r_t > 0$. Then~(\ref{skI2}) says that $\theta_t$ equals $-\pi(2t+a-1)(k-1)/p$ modulo $2\pi$.
Increasing $t$ by $1$ rotates $\rho_t$ by $-2\pi(k-1)/p$. Using the fact that $k-1$ is invertible modulo $p$, we have the following.
If $(a-1)(k-1)$ is even, then the set of $\theta_t$ for $t \in \Z p$ is precisely $0,2\pi/p,\ldots,(2p-2)\pi/p$, so there is a unique $t$ (resp.\ a unique $t'$) in $\Z p$ for which $\theta_t=\pi+\pi/p$ (resp.\ $\theta_{t'} = \pi-\pi/p$). Furthermore, $t' = -(a-1)-t$ and
$I+t' = -(I+t)$; thus $I+t$ and $I+t'$ have the same set of dilations.
If $(a-1)(k-1)$ is odd, then the set of $\theta_t$ for $t \in \Z p$ is precisely $\pi/p,3\pi/p,\ldots,(2p-1)\pi/p$, so there is a unique $t \in \Z p$ for which $\theta_t = \pi$.
We call $t$ (and $t'$, if it exists) \emph{optimal}.

Let $t$ be optimal.
To prove the theorem, we will show that $F(\xi\cdot(I+t)) < F(A)$ (and so $s_k(\xi\cdot(I+t))<s_k(A)$) for  any $a$-subset $A\subseteq\Z p$ which is not a dilation of $I+t$.

We will first show that $F(I+t)<F(A)$ for any $a$-subset $A$ which is not affine equivalent to an interval.
By Lemma~\ref{lm:MaxNontriv}(i), we have that $|\fourier{\ind_{I+t}}(\pm 1)|=m_1(a)$ and $\rho(A) \leq m_2(a)$.
Let $m_2'(a)$ be the maximum of $\fourier{\ind_J}(\gamma)$ over all length-$a$ intervals $J$ and $\gamma \in [2,p-2]$. 
Lemma~\ref{lm:MaxNontriv}(i) implies that $m_2'(a)<m_1(a)$.
Thus
\begin{eqnarray}\label{knot1eq2}
\left|F(I+t)-2(m_1(a))^{k+1}\cos(\theta_t) - F(A)\right| \leq (p-1)(m_2(a))^{k+1} + (p-3)\left(m_2'(a)\right)^{k+1}.
\end{eqnarray}
Now $\cos(\theta_t) \leq \cos(\pi-\pi/p) < -0.9$ since $p \geq 7$. 
This together with the fact that $k\ge k_0(a,p)$ and Lemma~\ref{lm:MaxNontriv} imply that the absolute value of $2(m_1(a))^{k+1}\cos(\theta_t)<0$  is greater than the right-hand size of~(\ref{knot1eq2}).
Thus $F(I+t) < F(A)$, as required.

The remaining case is when $A=\zeta\cdot(I+v)$ for some non-optimal $v \in \Z p$ and non-zero $\zeta \in \Z p$. Since $s_k(A)=s_k(I+v)$, we may assume that $\zeta=1$.
Note that $\cos(\theta_t) \leq \cos(\pi-\pi/p) < \cos(\pi-2\pi/p) \leq \cos(\theta_v)$.
Thus
\begin{align*}
F(I+t)-F(I+v) &\leq 2(m_1(a))^{k+1}(\cos(\theta_t)-\cos(\theta_v)) + (2p-4)(m_2'(a))^{k+1}\\
&\leq 2(m_1(a))^{k+1}(\cos(\pi-\pi/p)-\cos(\pi-2\pi/p)) + (2p-4)(m_2'(a))^{k+1} <0
\end{align*}
where the last inequality uses the fact that $k$ is sufficiently large.
Thus $F(I+t)<F(I+v)$, as required.
\epf

Finally, using similar techniques, we prove Theorem~\ref{th:k1}.

\bpf[Proof of Theorem~\ref{th:k1}.] 
Recall that $p\ge 7$, $a\in [3,p-3]$ and  $k > k_0(a,p)$ is sufficiently large with $k\equiv1\pmod p$. Let $I:=[a]$ and $I'=[a-1]\cup\{a\}$.

Suppose first that $a$ and $k$ are both even. Let $A\subseteq\Z p$ be
an arbitrary $a$-set not affine equivalent to the interval~$I$.
By Lemma~\ref{lm:MaxNontriv}, $I$ attains $m_1(a)$ (exactly at $x=\pm 1$), while $\rho(A)<m_1(a)$. Also, $m_2'(a)<m_1(a)$, where $
 m_2'(a):=\max_{\gamma\in [2,p-2]}|\fourier{\ind_I}(\gamma)|$. Thus
 \begin{eqnarray*}
F(I) - F(A) &\stackrel{(\ref{eq:SkF1}),(\ref{skI})}{\leq}& 2\sum_{\gamma=1}^{\frac{p-1}{2}} (-1)^\gamma\left|\fourier{\ind_I}(\gamma)\right|^{k+1} + \sum_{\gamma=1}^{p-1} \left|\fourier{\ind_A}(\gamma)\right|^{k+1}\\
 &\leq& -2(m_1(a))^{k+1} + (2p-4) (\max\{m_2(a),m_2'(a)\})^{k+1}\ <\ 0,
 \end{eqnarray*}
 where the last inequality uses the fact that $k$ is sufficiently large.
So  $s_k(a)=s_k(I)$.
Using Lemma~\ref{lm:MaxNontriv}, the same argument shows that, for all $B \in \binom{\Z p}{a}$, we have $s_k(B)=s_k(a)$ if and only if $B$ is an affine image of $I$.
This completes the proof of Part 1 of the theorem.

Suppose now that at least one of $a,k$ is odd. Let $A$ be an $a$-set not equivalent to~$I$. Again by Lemma~\ref{lm:MaxNontriv}, we have
 \begin{eqnarray*}
F(I) - F(A) &\geq&  \sum_{\gamma=1}^{p-1} \left|\fourier{\ind_I}(\gamma)\right|^{k+1} - \sum_{\gamma=1}^{p-1} \left|\fourier{\ind_A}(\gamma)\right|^{k+1}\\ 
&\geq& 2(m_1(a))^{k+1} - (p-1)(m_2(a))^{k+1} \ >\ 0.
\end{eqnarray*}
So the interval $I$ and its affine images have in fact the largest number of additive $(k+1)$-tuples among all $a$-subsets of $\Z p$. In particular, $s_k(a) < s_k(I)$.

Suppose that there is some $A \in \binom{\Z p}{a}$ which is not affine equivalent to $I$ or~$I'$. (If there is no such $A$, then the unique extremal sets are affine images of $I'$ for all $k > k_0(a,p)$, giving the required.)
Write $\rho := re^{\theta i} = \fourier{\ind_{I'}}(1)$.
Then by Lemma~\ref{lm:MaxNontriv}(ii), we have $r=m_2(a)$, and $\rho(A) \leq m_3(a)$.
Given $k \geq 2$, let $s \in \mathbb{N}$ be such that $k=sp+1$. Then
\begin{equation}\label{FI'}
 \Big|F(I') - 2m_2(a)^{k+1}\cos (sp\theta)-F(A)\Big|\le (p-1)m_3(a)^{k+1}+(p-3)\left(m_2'(a)\right)^{k+1}.
\end{equation}
Proposition~\ref{angleprop} implies that there is an even integer $\ell \in \I N$
for which $c := p\theta - \ell\pi \in (-\pi,\pi)\setminus\{0\}$. Let $\varepsilon := \frac{1}{3}\min\lbrace |c|,\pi-|c|\rbrace > 0$.
Given an integer $t$, say that $s\in\I N$ is \emph{$t$-good} if $sc \in ((t-\frac{1}{2})\pi+\varepsilon,(t+\frac{1}{2})\pi-\varepsilon)$.
This real interval has length $\pi-2\varepsilon > |c|>0$, so must contain at least one integer multiple of $c$.
In other words, for all $t \in \mathbb{Z}\setminus\{0\}$ with the same sign as $c$, there exists a $t$-good integer $s> 0$.
As $sp\theta\equiv sc\pmod{2\pi}$, the sign of $\cos(sp\theta)$ is $(-1)^{t}$. Moreover, Lemma~\ref{lm:MaxNontriv} implies that
$m_2(a)> m_3(a), m'_2(a)$.
Thus, when $k=sp+1>k_0(a,p)$, the absolute value of $2m_2(a)^{k+1}\cos(sp\theta)$ is greater than the right-hand side of~(\ref{FI'}).
Thus, for large $|t|$, we have  $F(A) > F(I')$ if $t$ is even and $F(A)<F(I')$ if $t$ is odd, implying the theorem by~\eqref{eq:SkF1}.
\epf

\section{Acknowledgements}

We are grateful to an anonymous referee for their careful reading of the manuscript and for many helpful comments and to Vsevolod Lev for informing us about~\cite{Lev01duke}.

\end{document}